\documentclass[11pt,a4paper]{article}
\usepackage[utf8]{inputenc}
\usepackage{amsmath}
\usepackage{amsfonts}
\usepackage{amssymb}
\usepackage{bm}
\usepackage{pdflscape}
\usepackage[margin=1in,footskip=0.25in]{geometry}
\usepackage{graphicx}
\usepackage{multirow}
\usepackage{graphicx}
\usepackage{wrapfig}
\usepackage{setspace}
\usepackage{lscape}
\usepackage[table,xcdraw]{xcolor}
\usepackage{float}

\usepackage{longtable}
\usepackage{courier}
\usepackage{listings}
\usepackage[printwatermark]{xwatermark}
\begin{document}

\noindent\textbf{\Large{A new method for constructing continuous distributions on the unit interval}}\\

\begin{center}
Aniket Biswas\\
Department of Statistics\\
Dibrugarh University\\
Dibrugarh, Assam, India-786004\\
Email:\textit{biswasaniket44@gmail.com}
\end{center}

\begin{center}
Subrata Chakraborty\\
Department of Statistics\\
Dibrugarh University\\
Dibrugarh, Assam, India-786004\\
Email:\textit{subrata\_stats@dibru.ac.in}
\end{center}

\section*{Abstract}
A novel approach towards construction of absolutely continuous distributions over the unit interval is proposed. Considering two absolutely continuous random variables with positive support, this method conditions on their convolution to generate a new random variable in the unit interval. This approach is demonstrated using some popular choices of the positive random variables such as the exponential, Lindley, gamma. Some existing distributions like the uniform and the beta are formulated with this method. Several new structures of density functions having potential for future application in real life problems are also provided. One of the new distributions having one parameter is considered for parameter estimation and real life modelling application and shown to provide better fit than the popular one parameter Topp-Leone model.\\

\noindent\textbf{\large{Keywords:}} Topp-Leone, Uniform, Beta, Exponential, Gamma, Lindley, Convolution, Ferlie-Gumble- Morgenstern copula.\\

\noindent\textbf{\large{MSC 2010:}} 60E05. \\

\section{Introduction}
Consider a departmental store where customers arrive independently of each other. Suppose the store opens at 10:00 a.m. and the second customer arrives at 11:00 a.m. Arrival time of the first customer was not observed. Then the store manager may be interested in finding the probability of the first customer arriving within half an hour from opening the store. If the customers arrive according to a Poisson process, then the probability of the first customer arriving within 10:30 a.m. is $0.5$. In fact, for given arrival time of the second customer, arrival time of the first customer is conditionally uniform over the duration that is uniform over 10:00 a.m to 11:00 a.m. This follows from the assumption that the inter-arrival times of the customers are independently and identically distributed (iid) exponential random variates. However, the conditional distribution does not remain uniform when any one of the above  assumptions is relaxed. Still the resulting conditional distribution has support $(0,1)$, the unit interval. This simple yet previously unaddressed approach clearly opens up a scope for constructing new distributions in the unit interval.\\ 

Importance of distributions with support $(0,1)$ is well established for modelling proportions, scores, rates, indices etc. Responses in the unit interval are often focus of investigation in many fields of application including finance (Gómez-Déniz et al. 2014, Biswas et al. 2020), public health (Mazucheli et al. 2019, Biswas and Chakraborty 2019), demography (Andreopoulos et al. 2019). The first model that comes to the mind of practitioners is obviously the well known Beta distribution (Johnson et al. 1995). Besides this, the Kumaraswamy's distribution (Kumaraswamy 1980, Jones 2009) and the Topp-Leone distribution (Topp and Leone 1955) also received attention. In recent times, a number of new models have been proposed as alternatives to the Beta and Kumaraswamy's distributions to enhance modelling flexibility as par context (Grassia 1977, Tadikamalla and Johnson 1982, Barndorff-Nielsen and Jorgensen 1991, Lemonte et al. 2013, Pourdarvish et al. 2015, Mazucheli et al. 2018, 2019, 2020). This recent surge in the number of research papers devoted to proposing new distribution in the unit interval clearly exhibits their growing relevance.\\

Most of the above distributions are generated through a suitable variable transformation of a baseline distribution with positive support. For $X$ being a positive valued random variable, the transformations $U=X/(1+X)$, $U=1/(1+X)$, $U=\exp(-X)$ are popularly used to derive a distribution with support $(0,1)$. Besides this, a general approach based on the cumulative distribution function-quantile function by Smithson and Shou (2017), Rodrigues et al. (2019) and Cancho et al. (2020). The objective of this paper is to describe a general procedure for construction of distributions with support $(0,1)$ from conditional convolution  approach.\\

Rest of the article is organized as follows. The new method of deriving distributions on $(0,1)$ is presented in the next section. Section 3 is devoted for illustrating the proposed methodology. Section 4 provides some detail on a new one parameter distribution along with real life application. This article is concluded with a potential generalization and further scopes. 
\section{Method}
Consider a two component absolutely continuous random vector $(X,Y)$ with joint cumulative distribution function (cdf) $F_{X,Y}$ and joint probability density function (pdf) $f_{X,Y}$. Let $F_X$ and $f_X$ are the marginal cdf and the marginal pdf of $X$, respectively. Similarly, $F_Y$ and $f_Y$ denote the marginal cdf and the marginal pdf of $Y$. The conditional cdf and pdf of $X$ given the realization on $Y$ are $F_{X|Y}$ and $f_{X|Y}$, respectively. \\

We assume that both $X$ and $Y$ have positive support, that is $P(X>0)=P(Y>0)=1$. For $Z=X+Y$, the cdf of $Z$ can be expressed as

\begin{equation}\label{eq:Fz}
    F_Z(z)=\int_0^z F_{X|Y}(z-y) f_Y(y)\,dy.
\end{equation}
The corresponding pdf is 

\begin{equation}\label{eq:fz}
    f_Z(z)=\int_0^z f_{X,Y}(z-y,y)\,dy=\int_0^z f_{X|Y}(z-y) f_Y(y)\, dy.
\end{equation}
In particular, when $X$ and $Y$ are independent, the cdf and the pdf of $Z$ can be obtained by replacing $F_{X|Y}(z-y)$ by $F_X(z-y)$ in (\ref{eq:Fz}) and $f_{X|Y}(z-y)$ by $f_X(z-y)$ in (\ref{eq:fz}), respectively. \\

We are interested in the conditional distribution of $X$ given $Z=1$. The corresponding random variable be $U$ with support $(0,1)$, the unit interval. Now $P(U\leq u)=P(X\leq u|Z=1)=P(X\leq u , X+Y=1)/f_Z(1)$. $P(X\leq u , X+Y=1)$ can be seen as $\int_0^u f_{X,Y}(v,1-v)\, dv$ and hence the pdf of $U$ is   

\begin{equation}\label{eq:fu}
    f_U(u)=\frac{f_{X,Y}(u,1-u)}{f_Z(1)}, \,\, \textrm{for}\,\, 0<u<1.
\end{equation}
If $X$ and $Y$ are independently distributed, then $f_{X,Y}(u,1-u)$ can be replaced by $f_X(u)f_Y(1-u)$ in (\ref{eq:fu}).

It may be noted that, bivariate distributions on $(0,\infty)\times (0,\infty)$ are seldom naturally available. The common techniques to construct bivariate distributions from given marginals involve marginal transformation method, trivariate reduction and bivariate copula (Balakrishnan and Lai 2009). To demonstrate the proposed method for constructing distributions on the unit interval, some bivariate distributions will be considered in Section 3. The bivariate density we use here is derived using Ferlie-Gumble-Morgenstern (FGM) copula (Gumble 1960) when the marginal densities $f_X$ and $f_Y$ are available. The joint density function of $(X,Y)$ with additional parameter $\alpha\in [-1,1]$ using FGM copula is

\begin{equation}\label{eq:fgm}
    f_{X,Y}(x,y)=f_X(x)f_Y(y)[1+\alpha(2F_X(x)-1)(2F_Y(y)-1)].
\end{equation}

\section{Main results}
In this section we apply the proposed method to first derive distributions in $(0, 1)$ which are known followed by some new ones which to the best of our knowledge are not yet reported in the literature.

\subsection{Known distributions}
As outlined in Section 1, let $X$ and $Y$ follow exponential distribution independently with same parameter $\theta$. Then $f_X(x|\theta)=\theta\exp(-\theta x)$ and $f_Y(y|\theta)=\theta\exp(-\theta y)$. Putting the the expressions of $f_X(x|\theta)$ and $f_Y(y|\theta)$ in (\ref{eq:fz}), we get $f_Z(1|\theta)=\theta^2\exp(-\theta)$. Now, putting the expressions of $f_X(x|\theta)$, $f_Y(y|\theta)$ and $f_Z(1|\theta)$ in (\ref{eq:fu}) we get $f_U(u)=1$ for $0<u<1$ which is pdf of the uniform uniform distribution over unit interval. \\

Let $G(\alpha,\beta)$ denote the gamma distribution with density function $\alpha^\beta x^{\beta-1}\exp(-\alpha x)/\Gamma(\beta)$ for $\alpha,\beta >0$. Consider $X\sim$G$(\alpha,\beta_1)$ and $Y\sim$G$(\alpha,\beta_2)$ independently. Using (\ref{eq:fz}) it is easy to see that $Z\sim G(\alpha,\beta_1+\beta_2)$. Now, putting $f_X(x|\alpha,\beta_1)$, $f_Y(y|\alpha,\beta_2)$ and $f_Z(1)=\alpha^{\beta_1+\beta_2}\exp(-\alpha)/\Gamma(\beta_1+\beta_2)$ in (\ref{eq:fu}) we get $f_U(u|\beta_1,\beta_2)=\Gamma(\beta_1+\beta_2)u^{\beta_1-1}(1-u)^{\beta_2-1}/\Gamma(\beta_1)\Gamma(\beta_2)$ which is the pdf of beta distribution with parameters $\beta_1$ and $\beta_2$.\\

While the uniform and beta cases are easily seen, generating the other distributions in the unit interval as mentioned in Section 1 using the proposed method remain a challenge. However, we do not claim that the distributions can (or cannot) be generated using the proposed method.  

\subsection{New distributions}
Here we consider different distributions for $(X, Y)$ to construct a number of new distributions on the unit interval.\\

\noindent\textbf{Distribution 1:} Let $X$ and $Y$ follow exponential distribution independently with parameter $\theta_1$ and $\theta_2$, respectively. Putting $f_X(x|\theta_1)=\theta_1\exp(-\theta_1x)$ and $f_Y(y|\theta_2)=\theta_2\exp(-\theta_2y)$ in (\ref{eq:fz}) we get 
\begin{equation}
    f_Z(1)=\begin{cases}
    \frac{\theta_1\theta_2(e^{-\theta_1}-e^{-\theta_2})}{\theta_1-\theta_2}\quad \textrm{for}\quad \theta_1\neq\theta_2\\
    \theta^2e^{-\theta} \quad \textrm{for}\quad \theta_1=\theta_2=\theta.
    \end{cases}
\end{equation}
Now, putting $f_X(x|\theta_1)$, $f_Y(y|\theta_2)$ and $f_Z(1)$ in (\ref{eq:fu}) and replacing $\theta_1-\theta_2$ by $\delta$ we obtain 
\begin{equation} \label{eq:d1}
f_U(u)=\begin{cases}
\frac{\delta e^{-\delta u}}{1-e^{-\delta}}\quad \textrm{for} \quad \delta\neq 0\\
1\quad \textrm{for} \quad \delta=0.
\end{cases}
\end{equation}
The density function in (\ref{eq:d1}) reduces to that of truncated exponential distribution with region of truncation $(1,\infty)$ when we restrict the parameter space to $\delta>0$. However, the new genesis does not require the parameter $\delta$ to be positive. Thus the density in (\ref{eq:d1}) provides more flexibility with $\delta\in\mathcal{R}$, the real line despite having the same mathematical structure as the said truncated exponential distribution. From Figure 1 it is clear that the density can be both increasing and decreasing depending on $\delta$. For $\delta<0$, it is increasing with mode at $1$ and for $\delta>0$, it is decreasing with mode at $0$.\\

\begin{figure}[H]
\begin{center}
\includegraphics[height=3in,width=4in,angle=0]{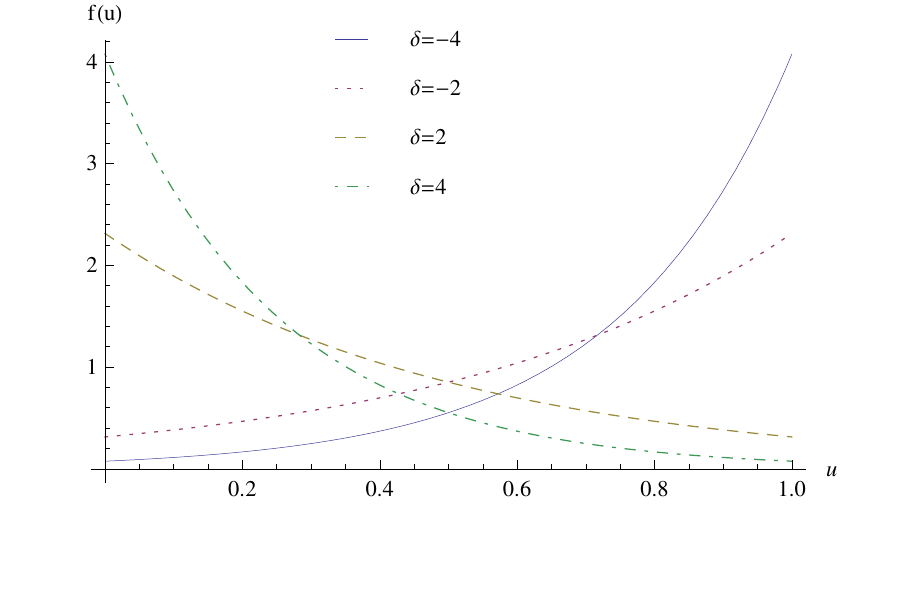}
\caption{Plot of the density function of Distribution 1.}
\end{center}
\end{figure}

\noindent\textbf{Distribution 2:} Let $X$ and $Y$ both marginally follow exponential distribution with parameters $\theta$. Putting $f_X(x|\theta)=\theta\exp(-\theta x)$, $F_X(x|\theta)=1-\exp(-\theta x)$, $f_Y(y|\theta)=\theta\exp(-\theta y)$ and $F_Y(y|\theta)=1-\exp(-\theta y)$ in (\ref{eq:fgm}) we obtain the following bivariate exponential distribution with parameters $\theta>0$ and $-1\leq\alpha\leq 1$.
\begin{equation}\label{eq:fgmd2}
    f_{X,Y}(x,y|\theta,\alpha)=\theta^2 e^{-\theta(x+y)}\left[1+\alpha(1-2e^{-\theta x})(1-2e^{-\theta y})\right]
\end{equation}
Using (\ref{eq:fgmd2}) in (\ref{eq:fz}) we obtain
\begin{equation}
    f_Z(1)=\theta^2 e^{-\theta}\left[1+\alpha\left\{1+4e^{-\theta}-\frac{4}{\theta}+4\frac{e^{-\theta}}{\theta}\right\}\right].
\end{equation}
Putting (\ref{eq:fgmd2}) and (\ref{eq:fz}) in (\ref{eq:fu}) we obtain the following density function.
\begin{equation}\label{eq:d2}
    f_U(u)=\frac{1+\alpha(1-2e^{-\theta u})(1-2e^{-\theta(1-u)})}{\left[1+\alpha\left\{1+4e^{-\theta}-\frac{4}{\theta}+4\frac{e^{-\theta}}{\theta}\right\}\right]}
\end{equation}
This distribution has two parameters and may prove to be useful in modelling applications. It is clear that for $\alpha=0$, $X$ and $Y$ are independent and consequently (\ref{eq:d2}) reduces to uniform distribution. From Figure 2, it can be seen that the distribution is symmetric about $0.5$ and it possesses both bath-tub shaped and bell-shaped curves.\\

\begin{figure}[H]
\begin{center}
\includegraphics[height=3in,width=6in,angle=0]{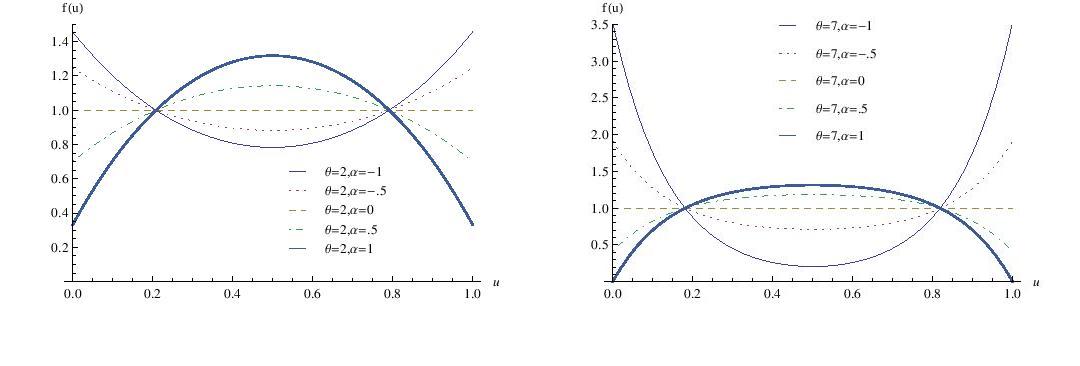}
\caption{Plot of the density function of Distribution 2.}
\end{center}
\end{figure}

\noindent\textbf{Distribution 3:} Let $X$ and $Y$ follow Lindley distribution independently and identically with parameter $\theta>0$. Putting $f_X(x|\theta)=\theta^2(1+x)\exp(-\theta x)/(1+\theta)$ and $f_Y(y|\theta)=\theta^2(1+y)\exp(-\theta y)/(1+\theta)$ in (\ref{eq:fz}) we get
\begin{equation}\label{eq:LLZ}
    f_Z(1)=\frac{13}{6}\frac{e^{-\theta}\theta^4}{(1+\theta)^2}.
\end{equation}
As in previous cases, using the density functions of $X$ and $Y$ along with (\ref{eq:LLZ}) in (\ref{eq:fu}) we obtain
\begin{equation}\label{eq:LLU}
    f_U(u)=\frac{6}{13}(2-u)(1+u).
\end{equation}
The resulting distribution is free from any parameter and thus not flexible as clear from Figure 3. Clearly (\ref{eq:LLU}) is symmetric about $0.5$.\\

\begin{figure}[H]
\begin{center}
\includegraphics[height=3in,width=4in,angle=0]{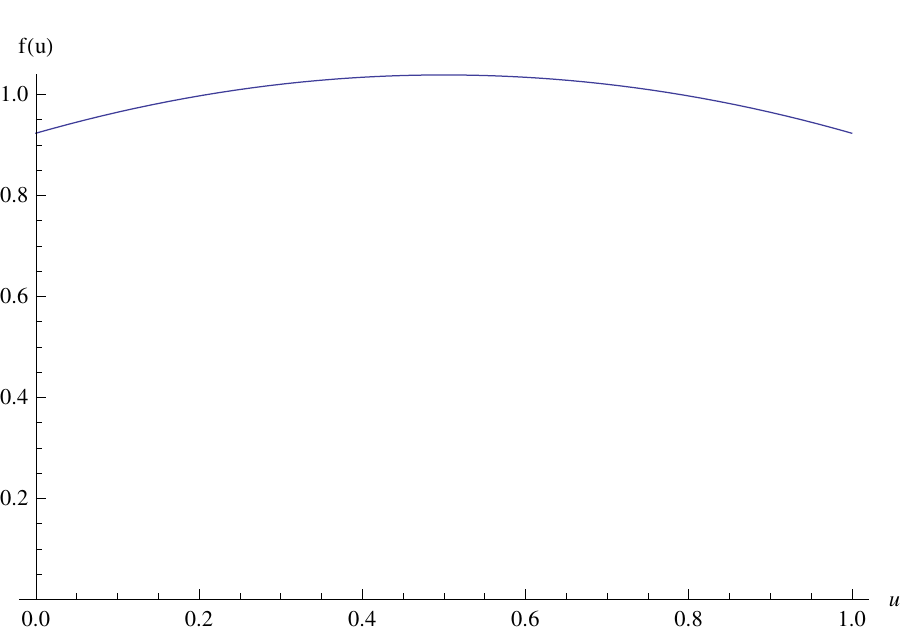}
\caption{Plot of the density function of Distribution 3.}
\end{center}
\end{figure}

\noindent\textbf{Distribution 4:} Let $X$ and $Y$ follow Lindley distribution independently with parameters $\theta_1$ and $\theta_2$, respectively. Here $f_X(x|\theta_1)=\theta_1^2(1+x)\exp(-\theta_1 x)/(1+\theta_1)$ and $f_Y(y|\theta_2)=\theta_2^2(1+y)\exp(-\theta_2 y)/(1+\theta_2)$. The density function of $Z=X+Y$ is evaluated using (\ref{eq:fz}) but we do not report it here due to its messy structure. However, the density function of $U$ obtained using $f_X(x|\theta_1)$, $f_Y(y|\theta_2)$ and $f_Z(1)$ in (\ref{eq:fu}) is 
\begin{equation}\label{eq:LLU2}
    f_U(u)=\begin{cases}
    \frac{\delta^3 e^{-\delta u} (u-2)(1+u)}{(2-\delta-2\delta^2)+e^{-\delta}(2\delta^2-2-\delta)}\quad \textrm{for} \quad \delta\neq0\\
    \frac{6}{13}(2-u)(1+u)\quad \textrm{for} \quad \delta=0
    \end{cases}
\end{equation}
for $\delta=\theta_1-\theta_2\in \mathcal{R}$. Obviously the density in (\ref{eq:LLU2}) having one parameter is an extension of (\ref{eq:LLU}). From Figure 4 and Figure 1, it can be seen that the Distribution 4 and Distribution 1 behave quite similarly.\\

\begin{figure}[H]
\begin{center}
\includegraphics[height=3in,width=4in,angle=0]{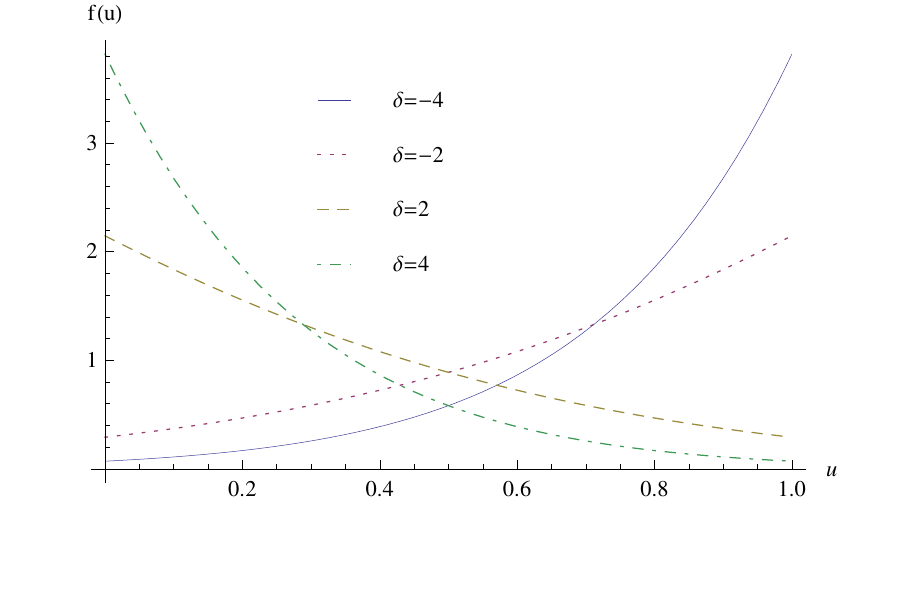}
\caption{Plot of the density function of Distribution 4.}
\end{center}
\end{figure}

\noindent\textbf{Distribution 5:} Let $X$ follow exponential distribution with parameter $\theta$ and $Y$ follow Lindley distribution with parameter $\theta$ independently. Here $f_X(x|\theta)=\theta\exp(-\theta x)$ and $f_Y(y|\theta)=\theta^2(1+y)\exp(-\theta y)/(1+\theta)$. From (\ref{eq:fu}) we get
\begin{equation}\label{eq:EL1}
    f_U(u)=\frac{2}{3}(2-u).
\end{equation}
This distribution has no parameter and thus not flexible for modelling purposes . The density in (\ref{eq:EL1}) is decreasing and it decreases from $4/3$ at $u=0$ to $2/3$ at $u=1$ linearly as given in Figure 5.\\

\begin{figure}[H]
\begin{center}
\includegraphics[height=3in,width=4in,angle=0]{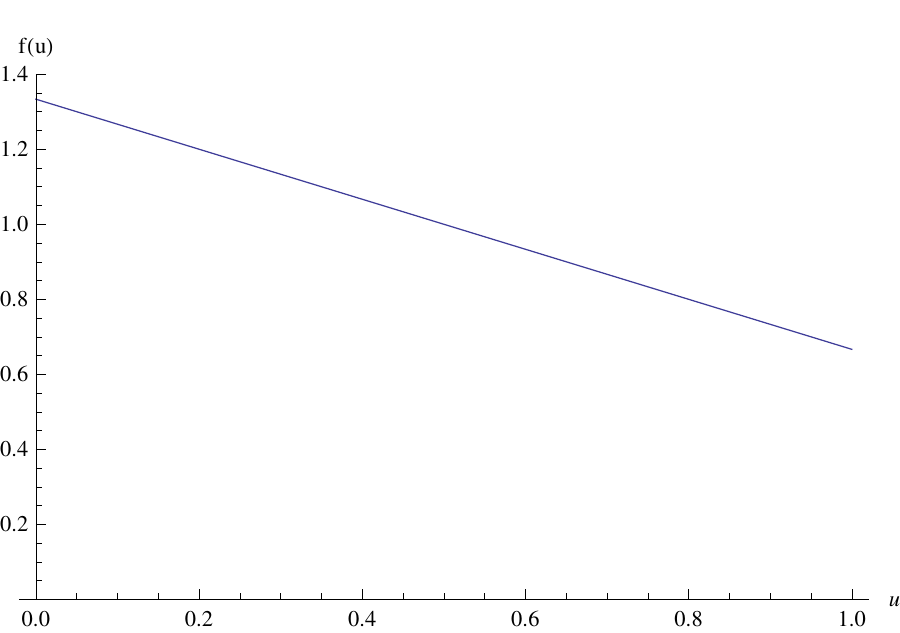}
\caption{Plot of the density function of Distribution 5.}
\end{center}
\end{figure}

\noindent \textbf{Distribution 6:} Let $X$ follow exponential distribution with parameter $\theta_1$ and $Y$ follow Lindley distribution with parameter $\theta_2$ independently. Here $f_X(x|\theta_1)=\theta_1\exp(-\theta_1 x)$ and $f_Y(y|\theta_2)=\theta_2^2(1+y)\exp(-\theta_2 y)/(1+\theta_2)$. Using the density functions first in (\ref{eq:fz}) and then in (\ref{eq:fu}) we get
\begin{equation}\label{EL2}
    f_U(u)=\begin{cases}
    \frac{\delta^2 e^{\delta(1-u)}(2-u)}{1-\delta+(2\delta-1)e^\delta}\quad \textrm{for}\quad \delta\neq 0\\
    \frac{2}{3}(2-u)\quad \textrm{for}\quad \delta= 0.
    \end{cases}
\end{equation}
Here $\delta=\theta_1-\theta_2\in \mathcal{R}$. A few representative plots of the density function are provided in Figure 6.\\

\begin{figure}[h]
\begin{center}
\includegraphics[height=3in,width=4in,angle=0]{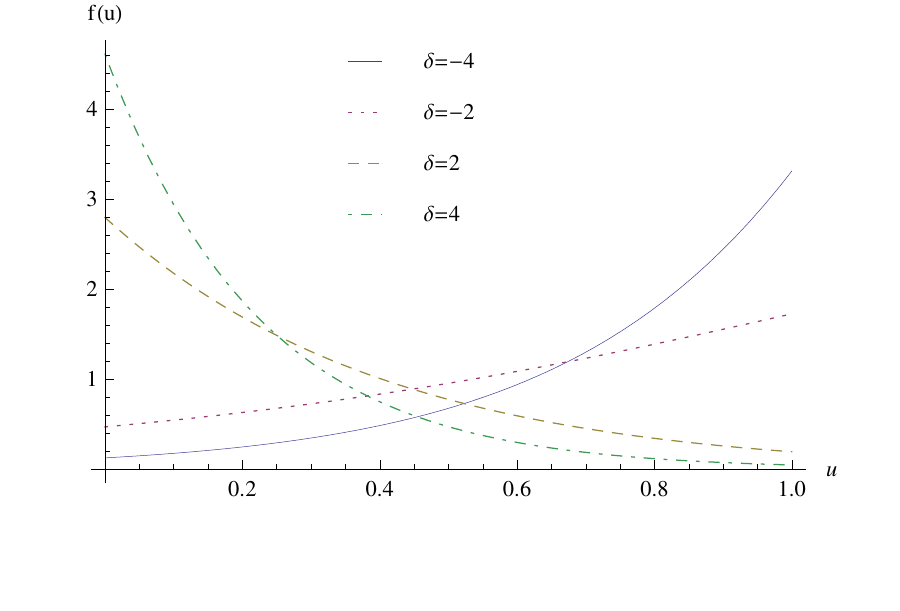}
\caption{Plot of the density function of Distribution 6.}
\end{center}
\end{figure}

\noindent \textbf{Distribution 7:} Consider $X$ following exponential distribution with parameter $\theta$ and $Y$ following gamma distribution with parameters $\alpha$ and $\beta$. Putting $f_X(x|\theta)=\theta\exp(-\theta x)$ and $f_Y(y|\alpha,\beta)=\alpha^\beta y^{\beta-1}\exp(-\alpha x)/\Gamma(\beta)$ in (\ref{eq:fz}) we get
\begin{equation}\label{eq:egz}
    f_Z(1)=\begin{cases}
    \frac{\theta \alpha^\beta e^{-\theta}\{\Gamma(\beta)-\Gamma(\alpha-\theta)\}}{\Gamma(\beta) (\alpha-\theta)^\beta}\quad \textrm{for} \quad \alpha\neq \theta\\
    \frac{\alpha^{\beta+1}\exp(-\alpha)}{\Gamma(\beta+1)} \quad \textrm{for}\quad \alpha=\theta.
    \end{cases}
\end{equation}
Here $\Gamma(s,t)$ denotes the incomplete gamma function $\int_t^\infty v^{s-1} \exp(-v)\, dv$. Now using the density functions of $X$ and $Y$ along with (\ref{eq:egz}) in (\ref{eq:fu}) we obtain 
\begin{equation}\label{eq:egu}
    f_U(u)=\begin{cases}
    \frac{\delta^\beta}{\Gamma(\beta)-\Gamma(\beta,\delta)}e^{\delta(1-u)}(1-u)^{\beta-1}\quad \textrm{for}\quad \beta>0,\delta\neq 0\\
    \beta (1-u)^{\beta-1}\quad \textrm{for} \quad \alpha>0,\delta=0
    \end{cases}
\end{equation}
where $\delta=\alpha-\theta\in \mathcal{R}$. The density in (\ref{eq:egu}) has two parameters and may prove to be a competitor of the popular beta distribution as evident from Figure 7.\\

\begin{figure}[H]
\begin{center}
\includegraphics[height=2in,width=6in,angle=0]{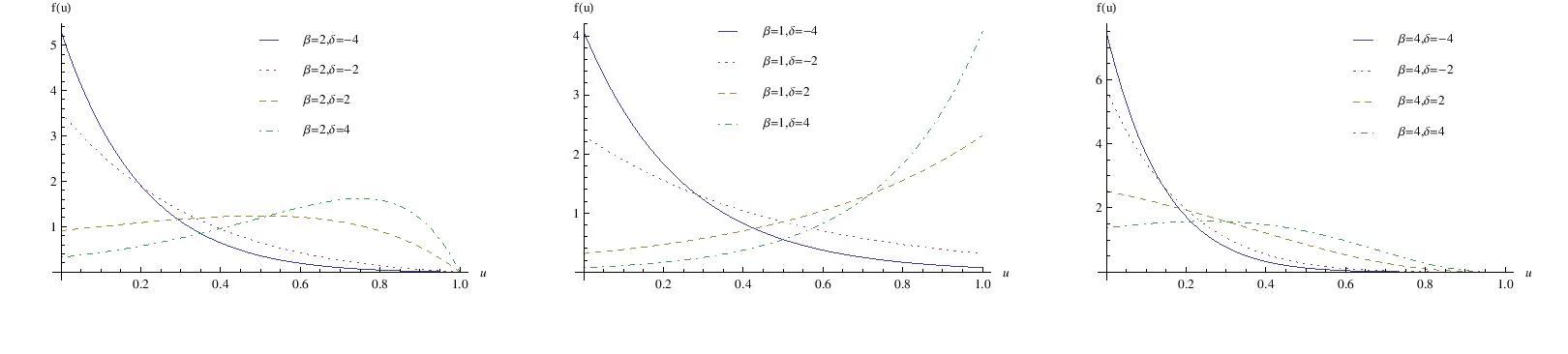}
\caption{Plot of the density function of Distribution 7.}
\end{center}
\end{figure}

\noindent \textbf{Distribution 8:} Let $X$ be a Lindley random variable with parameter $\theta$ and $Y$ be a gamma random variable with parameters $\theta$ and $\beta$. Here $f_X(x|\theta)=\theta^2(1+x)\exp(-\theta x)/(1+\theta)$ and $f_Y(y|\theta,\beta)=\theta^\beta y^{\beta-1}\exp(-\theta x)/\Gamma(\beta)$. Putting the densities of $x$ and $Y$ first in (\ref{eq:fz}) and then using $f_Z(1)$ along with the density functions in (\ref{eq:fu}) we obtain
\begin{equation}\label{eq:flg}
    f_U(u)=\frac{\beta(1+\beta)}{2+\beta} (1+u) (1-u)^{\beta-1}\quad \textrm{for} \quad \beta>0.
\end{equation}
The parameter $\theta$ labelling the distributions of both $X$ and $Y$ surprisingly plays no role in (\ref{eq:flg}). Despite being a one-parameter distribution, it is quite flexible for modelling purposes since the corresponding density function has a wide range of shapes as can be seen from Figure 8. The derived distribution can also be seen as a mixture of two beta distributions where the mixing proportions are also functions of the parameter $\beta$ .

\begin{figure}[H]
\begin{center}
\includegraphics[height=3in,width=4in,angle=0]{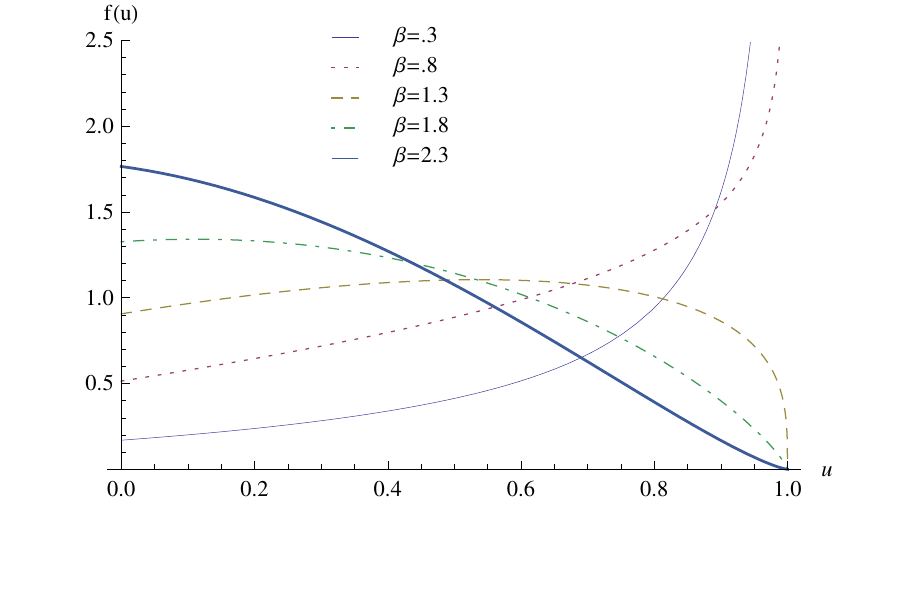}
\caption{Plot of the density function of Distribution 8.}
\end{center}
\end{figure} 

\section{LCG: A new one parameter distribution}
Distribution 8 is derived from the distribution of Lindley random variable conditioned on the convolution of the same with gamma random variable. Thus a working abbreviation for the distribution with pdf in (\ref{eq:flg}) is taken as $LCG(\beta)$ and we write $U\sim LCG(\beta)$. The cumulative distribution function (cdf) of $LCG(\beta)$ is
\begin{equation}\label{eq:Flg}
    F_U(u)=1-\frac{(1-u)^\beta(2+\beta+u\beta)}{2+\beta},\quad  0<u<1.
\end{equation}
Figure 9 exhibits the functional form of $F_U$ for different choices of $\beta$. With increasing $\beta$, the cdf becomes more and more convex in shape. 
\begin{figure}[H]
\begin{center}
\includegraphics[height=3in,width=4in,angle=0]{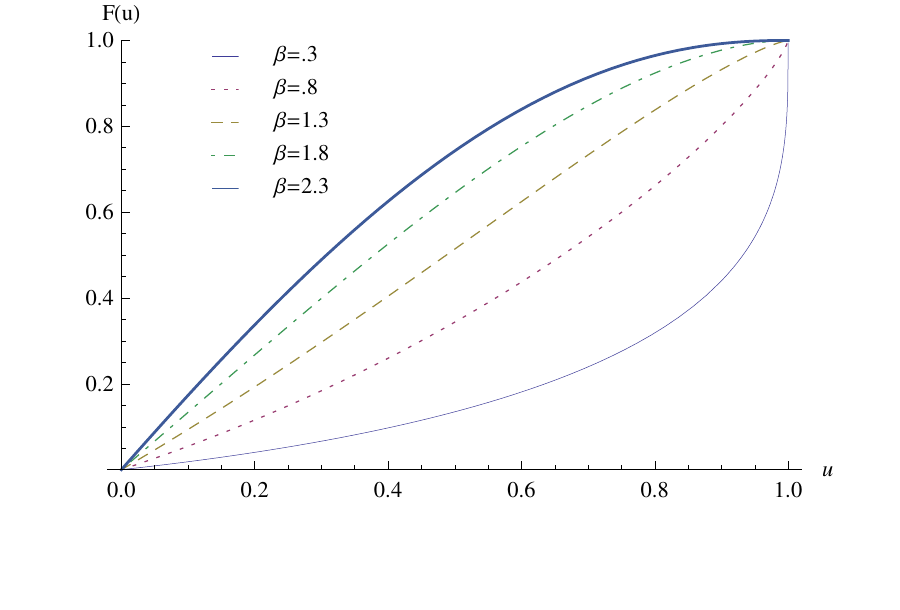}
\caption{Plot of the cumulative distribution function of LCG distribution.}
\end{center}
\end{figure}
The corresponding hazard rate is 
\begin{equation}\label{eq:hlg}
    h(u)=\frac{\beta(1+\beta)}{2+\beta+u\beta}\frac{1+u}{1-u},\quad   0<u<1.
\end{equation}
The LCG family has increasing failure rate as shown in Figure 10.
\begin{figure}[H]
\begin{center}
\includegraphics[height=3in,width=4in,angle=0]{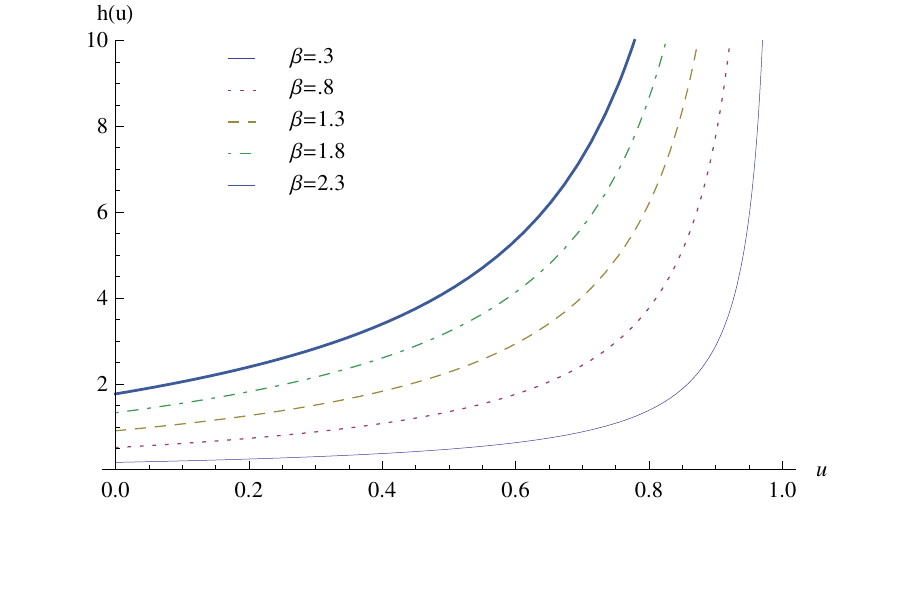}
\caption{Plot of the hazard rate of LCG distribution.}
\end{center}
\end{figure}
\noindent The lower order moments are very simple, the mean and the variance of $LCG(\beta)$ is given below.
\begin{equation}
    E(U)=\frac{4+\beta}{(2+\beta)^2}\quad \textrm{and} \quad V(U)=\frac{\beta(16+9\beta+\beta^2)}{(2+\beta)^4(3+\beta)}
\end{equation}

\subsection{Generating random observations}
To generate random observations from $LCG(\beta)$ it is convenient to use the cdf in (\ref{eq:Flg}). We employ cdf inversion technique as the structure is not too complex. For a random observation $v$ from uniform distribution over $(0,1)$ we equate $F_U(u)=v$ and find solution of the equation in $u$. This $u$ is a random observation from $LCG(\beta)$. This technique produces the following nonlinear equation in $u$
\begin{equation}\label{eq:Flgv}
    (1-u)^\beta(2+\beta+u\beta)-(1-v)(2+\beta)=0.
\end{equation}
We solve (\ref{eq:Flgv}) by \textit{uniroot} function available on basic \textbf{R}.
\subsection{Estimation}
Since $E(U)$ is structurally simple, a method of moment estimator can easily be implemented. For a given dataset $\bm{u}=(u_1, u_2, ..., u_n)$, we calculate the sample mean $\Bar{u}=(1/n)\sum_{i=1}^n u_i$ and form the following equation.
\begin{equation}\label{eq:mme1}
    E(U)=\frac{4+\beta}{(2+\beta)^2}=\Bar{u}
\end{equation}
The equation in (\ref{eq:mme1}) can be rewritten as 
\begin{equation}\label{eq:mme2}
    \Bar{u}\beta^2+(4\Bar{u}-1)\beta+4(\Bar{u}-1)=0.
\end{equation}
The method of moment estimator $\hat{\beta}_{MM}$ is the solution of (\ref{eq:mme2}) in $\beta$. \\

The log-likelihood function of $\beta$ given the dataset $\bm{u}$ is 
\begin{equation}\label{mle1}
    l(\beta|\bm{u})=n\log\left[\frac{\beta(1+\beta)}{2+\beta}\right]+\sum_{i=1}^n \log(1+u_i)+(\beta-1)\sum_{i=1}^n \log(1-u_i).
\end{equation}
Differentiating $l(\beta|\bm{u})$ with respect to $\beta$ and equating it with $0$ yields the following.
\begin{equation}\label{eq:mle2}
    T(\bm{u})\beta^2+(3T(\bm{u})-1)\beta+(2T(\bm{u})-4)=0
\end{equation}
Here $T(\bm{u})=-(1/n)\sum_{i=1}^n \log (1-u_i)$.  To solve both (\ref{eq:mme2}) and (\ref{eq:mle2}) in $\beta$, we use \textit{polyroot} function available in basic \textbf{R}. The solution of (\ref{eq:mle2}) is the maximum likelihood estimate namely $\hat{\beta}_{ML}$.
\begin{center}
\begin{table}[h!]
\centering
\vspace{1cm}
	\caption{Bias and mean-squared error (MSE) of $\hat{\beta}_{MM}$ and $\hat{\beta}_{ML}$.}
	\label{tab:simul}
	\begin{tabular}{cccccc}
		\hline\\
		$n$ & $\beta$ & Bias($\hat{\beta}_{MM}$) & MSE($\hat{\beta}_{MM}$) & Bias($\hat{\beta}_{ML}$) & MSE($\hat{\beta}_{ML}$)  \\ 
		\\\hline
		\multirow{3}{*}{20} & 2.30 & 0.07792 &0.27779  & -0.21959 &0.33440   \\
		& 7.80 &0.34196   &3.74806  &0.21788  &3.72484\\
		&  15.00&0.52377  &12.36922  &0.46433  &12.36009  \\ \hline
		\multirow{3}{*}{40} & 2.30 &0.01869  &0.12236  &-0.29453  &0.21342   \\
		& 7.80 &0.13403   &1.48331  &-0.00307  &1.48377\\
		&15.00 &0.38477   &5.71176  &0.30913  &5.72911  \\ \hline
				\multirow{3}{*}{60} &2.30  &0.01719  &0.08165  &-0.29909  &0.17421  \\
		& 7.80 &0.12145  &0.96437  &-0.02729  &0.96883 \\
		&15.00 & 0.16190  &3.51346  &0.08157&3.51912 \\ \hline
		
				\multirow{3}{*}{80} & 2.30 &0.00449 &0.06157 &-0.31263 &0.16230  \\
		&  7.80 &0.12302 &0.70574 &-0.03383 &0.70355  \\
		& 15.00 & 0.14265 & 2.69120 &0.05492  &2.69554  \\ \hline
		
				\multirow{3}{*}{100} & 2.30  & 0.01120 &0.04923  &-0.30954  &0.14639   \\
		& 7.80 &0.06890   &0.55969  &-0.08804  &0.57761\\
		& 15.00 & 0.17952 &2.35752  &0.08494  &2.35630  \\ \hline
		 
	\end{tabular}%
\end{table}
\end{center}
A simulation experiment is performed to study the relative performance of the estimators. We generate 1000 samples of different sizes ($n$) using the method in Section 4.1 and compute the average bias and mean squared error of the estimators. From Table 1, it is seen that the performance of both the estimators are quite similar. With increasing $\beta$, the mean squared error of both the estimators increase. As desired the mean squared error and absolute bias of both the estimators decrease with increase in sample size.
\subsection{Application}
Consider the following datasets from Caramanis et al. (1983) providing computation time of two different algorithms namely SC16 and P3. \\

\noindent SC16:= (0.853, 0.759, 0.866, 0.809, 0.717, 0.544, 0.492, 0.403, 0.344, 0.213, 0.116,
0.116, 0.092, 0.070, 0.059, 0.048, 0.036, 0.029, 0.021, 0.014, 0.011, 0.008, 0.006)\\

\noindent P3:= (0.853, 0.759, 0.874, 0.800, 0.716, 0.557, 0.503, 0.399, 0.334, 0.207, 0.118,
0.118, 0.097, 0.078, 0.067, 0.056, 0.044, 0.036, 0.026, .019, 0.014, 0.010)\\

\noindent Mazumdar and Gaver (1984) comparison of the mentioned algorithms for estimating the unit capacity factor. For a population level comparison one may refer to Genc (2013) where the authors fitted Topp-Leone model to both SC16 and P3 datasets for further analysis. It is to be noted that to improve the population level comparison, one need to select a better suited model. One parameter distributions on the unit interval is rare and the Topp-Leone distribution is a benchmark in this category. The introduced LCG distribution is a new addition to this category and hence a real life data driven comparison with Topp-Leone is sufficient to justify its utility. Topp-Leone model fits the datasets well as can be seen in Genc (2013). However, we compare fitting of the Topp-Leone model with that of the LCG in Table 2  and find that the proposed LCG model is far better as indicated by the corresponding AIC values. 
\begin{center}
\begin{table}[h!]
\centering
\vspace{1cm}
	\caption{Goodness of fit.}
	\label{tab:simul}
	\begin{tabular}{ccccc}
		\hline\\
		Dataset & Distribution & MLE & Log-likelihood & AIC  \\ 
		\\\hline
		\multirow{2}{*}{SC16}  & Topp-Leone & 0.5943 & -11.3660 & 25.7837   \\
		 & LCG & 1.9876 & 2.9459 & -3.8981 \\ \hline
		\multirow{2}{*}{P3} & Topp-Leone & 0.6778 & -10.9016 & 23.8032    \\
		
		  & LCG & 1.8646 &2.3009 &-2.6098\\ \hline
	\end{tabular}%
\end{table}
\end{center}

\section{Discussion}
A number of new distributions with simple yet flexible structures are derived in the current work. Throughout the article we discussed about the distribution of $U=X|X+Y=1$. It is obvious that in many situations the distribution of $V=Y|X+Y=1$ is different from $U$ and hence may be of interest fom new distribution perspective. However, we refrain from detailing such density function since the same can be done by simply finding density function of $V=1-U$. We envisage that the proposed method and the derived distributions will be a useful addition in the existing literature.  In fact the LCG is shown to be a better one parameter model than the benchmark Topp-Leone distribution. Since the LCG model possesses simple structure of the mean, modelling unit responses with covariates using LCG is of interest. Moreover Section 4.2 hints a productive study on the stress-strength reliability of the LCG distribution. These problems are under consideration and will be reported as independent works. Detailed investigation on the other derived distributions are required and further research on this is warranted. This article is intended to motivate researchers to delve into the possibilities of exploring new distributions with different choices of the joint distributions of $(X,Y)$. The proposed method can naturally be extended. One may consider the random vector $(X_1, X_2, ..., X_n)$ with $P(X_i>0)=$ for $i=1, 2, ..., n$. For $k<n$, $U=X_1+X_2+...+X_k|X_1+X_2+...+X_n=1$  is a random variable with support $(0,1)$. In our case $n=2$ and $k=1$.

\section*{References}
Andreopoulos, P., Bersimis, G. F., Tragaki, A., \& Rovolis, A. (2019). Mortality modeling using probability distributions. APPLICATION in greek mortality data. \textit{Communications in Statistics-Theory and Methods}, 48(1), 127-140.\\
\\
Balakrishnan, N., \& Lai, C. D. (2009). \textit{Continuous bivariate distributions}. Springer Science \& Business Media.\\
\\
Barndorff-Nielsen, O. E., \& Jørgensen, B. (1991). Some parametric models on the simplex. \textit{Journal of Multivariate Analysis}, 39(1), 106-116.\\
\\
Biswas, A., \& Chakraborty, S. (2019). $R= P (Y< X)$ for unit-Lindley distribution: inference with an application in public health. \textit{arXiv preprint arXiv:1904.06181}.\\
\\
Biswas, A., Chakraborty, S., \& Mukherjee, M. (2020). On estimation of stress–strength reliability with log-Lindley distribution. \textit{Journal of Statistical Computation and Simulation}, 1-23.\\
\\
Cancho, V. G., Bazán, J. L., \& Dey, D. K. (2020). A new class of regression model for a bounded response with application in the study of the incidence rate of colorectal cancer. \textit{Statistical Methods in Medical Research}, 29(7), 2015-2033.\\
\\
Caramanis, M., Stremel, J., Fleck, W., \& Daniel, S. (1983). Probabilistic production costing: an investigation of alternative algorithms. \textit{International Journal of Electrical Power \& Energy Systems}, 5(2), 75-86.\\
\\
Genc, A. I. (2013). Estimation of $P (X> Y)$ with Topp–Leone distribution. \textit{Journal of Statistical Computation and Simulation}, 83(2), 326-339.\\
\\
Gómez-Déniz, E., Sordo, M. A., \& Calderín-Ojeda, E. (2014). The Log–Lindley distribution as an alternative to the beta regression model with applications in insurance. \textit{Insurance: Mathematics and Economics}, 54, 49-57.\\
\\
Grassia, A. (1977). On a family of distributions with argument between 0 and 1 obtained by transformation of the gamma and derived compound distributions. \textit{Australian Journal of Statistics}, 19(2), 108-114.\\
\\
Gumbel, E. J. (1960). Multivariate distributions with given margins and analytical examples.\textit{ Bulletin de l’Institut International de Statistique}, 37(3), 363-373.\\
\\
Johnson, N. L., Kotz, S., \& Balakrishnan, N. (1995). \textit{Continuous Univariate Distributions}, volume 2. John Wiley\&Sons. Inc.,, 75.\\
\\
Jones, M. C. (2009). Kumaraswamy’s distribution: A beta-type distribution with some tractability advantages. \textit{Statistical Methodology}, 6(1), 70-81.\\
\\
Kumaraswamy, P. (1980). A generalized probability density function for double-bounded random processes. \textit{Journal of Hydrology}, 46(1-2), 79-88.\\
\\
Lemonte, A. J., Barreto-Souza, W., \& Cordeiro, G. M. (2013). The exponentiated Kumaraswamy distribution and its log-transform.\textit{ Brazilian Journal of Probability and Statistics}, 27(1), 31-53.\\
\\
Mazucheli, J., Menezes, A. F. B., \& Ghitany, M. E. (2018). The unit-Weibull distribution and associated inference. \textit{Journal of Applied Probability and Statistics}, 13, 1-22.\\
\\
Mazucheli, J., Menezes, A. F. B., Fernandes, L. B., de Oliveira, R. P., \& Ghitany, M. E. (2020). The unit-Weibull distribution as an alternative to the Kumaraswamy distribution for the modeling of quantiles conditional on covariates. \textit{Journal of Applied Statistics}, 47(6), 954-974.\\
\\
Mazumdar, M., \& D. P. Gaver. (1984). On the computation of power-generating system reliability indexes.\textit{ Technometrics} 26, no. 2 (1984): 173-185.\\
\\
Mazucheli, J., Menezes, A. F. B., \& Chakraborty, S. (2019). On the one parameter unit-Lindley distribution and its associated regression model for proportion data. \textit{Journal of Applied Statistics}, 46(4), 700-714.\\
\\
Pourdarvish, A., Mirmostafaee, S. M. T. K., \& Naderi, K. (2015). The exponentiated Topp-Leone distribution: Properties and application. \textit{Journal of Applied Environmental and Biological Sciences}, 5(7), 251-6.\\
\\
Rodrigues, J., Bazán, J. L., \& Suzuki, A. K. (2020). A flexible procedure for formulating probability distributions on the unit interval with applications. \textit{Communications in Statistics-Theory and Methods}, 49(3), 738-754.\\
\\
Smithson, M., \& Shou, Y. (2019). cdfquantreg: An R package for CDF-Quantile Regression.\\
\\
Tadikamalla, P. R., \& Johnson, N. L. (1982). Systems of frequency curves generated by transformations of logistic variables.\textit{ Biometrika}, 69(2), 461-465.\\
\\
Topp, C. W., \& Leone, F. C. (1955). A family of J-shaped frequency functions. \textit{Journal of the American Statistical Association}, 50(269), 209-219.

\section*{Appendix: R-codes}
\subsubsection*{C1: Random sample generation from LCG distribution.}
\begin{lstlisting}[language=R]
rGL=function(n,beta)
{
sample=0
for(i in 1:n)
{
v=runif(1)
f=function(u)
{
(1-u)^beta*(2+beta+u*beta)-(1-v)*(2+beta)
}
proc=uniroot(f,c(0,1))
sample[i]=proc$root
}
return(sample)
}
\end{lstlisting}

\subsubsection*{C2: Method of moment estimation for LCG distribution.}
\begin{lstlisting}[language=R]
MME=function(data)
{
ubar=mean(data)
z=c(4*(ubar-1),4*ubar-1,ubar)
return(max(Re(polyroot(z))))
}
\end{lstlisting}

\subsubsection*{C2: Maximum likelihood estimation for LCG distribution.}
\begin{lstlisting}[language=R]
MLE=function(data)
{
Tu=-mean(log(1-data))
z=c(2*Tu-4,3*Tu-1,Tu)
return(max(Re(polyroot(z))))
}
\end{lstlisting}
\end{document}